\documentclass[10pt,reqno]{amsart}
\usepackage{amssymb,amsmath,amscd}

\newcommand{\id}{\mathop{\rm id}}
\newcommand{\Innt}{\mathop{\rm Int}}
\newcommand{\loc}{\mathop{\rm loc}}
\newcommand{\vol}{\mathop{\rm vol}}
\newcommand{\const}{\mathop{\rm const}}

\newtheorem{theorem}{Theorem}

\newtheorem{lemma}{Lemma}

\begin{document}

\title[$L_{p,q}$-Cohomology of Warped Cylinders]
{$L_{p,q}$-Cohomology of Warped Cylinders}\thanks{Supported 
by INTAS (Grant 03--51--3251), the Specific Targeted
Project GALA within the NEST Activities of the Commission of the European
Communities (Contract No.~028766), the State Maintenance
Program for the Leading Scientific Schools
of the Russian Federation (Grant~NSh~8526.2006.1). }
\address{Yaroslav Kopylov
\newline\hphantom{iii} Sobolev Institute of Mathematics,
\newline\hphantom{iii} Akademik Koptyug Pr. 4,
\newline\hphantom{iii} 630090, Novosibirsk, Russia}%
\email{yakop@math.nsc.ru}
\maketitle
\begin{abstract}
We extend some results by Gol$'$dshtein, Kuz$'$minov, and Shvedov about
the $L_p$-cohomology of warped cylinders to $L_{p,q}$-cohomology for
$p\ne q$. As an application, we establish some sufficient conditions
for the nontriviality of the $L_{p,q}$-torsion of a surface of revolution in terms of some
Hardy constants.



\end{abstract}

\section{Introduction}
\label{intro}
Let $M$ be a Riemannian manifold. For $1\le p\le\infty$ and a positive continuous
function $\sigma:M\to\mathbb{R}$, denote by $L_p^j(M,\sigma)$ the Banach space of measurable
forms of degree~$j$ on~$M$ with the finite norm
$$
\|\omega\|_{L_p^j(M,\sigma)}=
\begin{cases}
 \biggl\{
 \int_M |\omega(x)|^p \sigma^p(x)\,dx \biggr\}^{1/p}
&\text{if $1\le p<\infty$,} \\
\operatornamewithlimits{ess\,sup\,}_{x\in M} |\omega(x)|\sigma(x)
&\text{if $p=\infty$.}
\end{cases}
$$
Here $dx$ stands for the volume form of $M$ and $|\omega(x)|$ is the modulus of the
exterior form $\omega(x)$. In the usual way, we also define the spaces $L_{p,\loc}(M)$.

Denote by $D^j(M)=C_0^{\infty,j}(M)$ the space of smooth forms of degree~$j$ on~$M$
having compact support included in~$\Innt M$.
A form
$\psi\in L_{1,\loc}^{j+1}(M)$
is called the ({\it weak}) {\it differential}
$d\omega$
of
$\omega\in L_{1,\loc}^j(M)$
if
$$
\int\limits_U \omega\land du
=(-1)^{j+1}\int\limits_U \psi\land u
$$
for every orientable domain
$U\subset\Innt M$
and every form
$u\in D^{\dim M-j-1}(M)$
having support in
$U$.

Put
$$
W_{p,q}^j(M,\sigma_j,\sigma_{j+1})=\{\omega\in L_p^j(M,\sigma_j)
\mid d\omega\in L_q^{j+1}(M,\sigma_{j+1})\}.
$$
The space $W_{p,q}^j(M,\sigma_j,\sigma_{j+1})$ is endowed with the norm
$$
\|\omega\|_{W_{p,q}^j(M,\sigma_j,\sigma_{j+1})}
=\|\omega\|_{L_p^j(M,\sigma_j)}+ \|d\omega\|_{L_q^{j+1}(M,\sigma_{j+1})}.
$$
If $p=q$ then it is often more convenient to consider the equivalent norm
$$
\|\omega\|_{W_p^j(M,\sigma_j,\sigma_{j+1})}
=\left(\|\omega\|^p_{L_p^j(M,\sigma_j)}+ \|d\omega\|^p_{L_p^{j+1}(M,\sigma_{j+1})}\right)^{1/p}.
$$
Denote by $V_{p,q}^j(M,\sigma_j,\sigma_{j+1})$
the closure of~$D^j(M)$ in the norm of $W_{p,q}^j(M,\sigma_j,\sigma_{j+1})$.

Given an arbitrary subset $A\subset M$, let $W_{p,q}^j(M,A,\sigma_j,\sigma_{j+1})$ be the
closure in $W_{p,q}^j(M,\sigma_j,\sigma_{j+1})$ of the subspace spanned by all forms
$\omega\in W_{p,q}^j(M,\sigma_j,\sigma_{j+1})$ which vanish on some neighborhood of~$A$
(depending on~$\omega$).

Let $Z_q^j(M,\sigma_j)$ be the subspace in~$W_{q,q}^j(M,\sigma_j)$ that consists of all forms
$\omega$ such that $d\omega=0$ and let
$$
B_{p,q}^j(M,\sigma_{j-1},\sigma_j)
= \{ \theta\in W_{q,q}^j(M,\sigma_j,\sigma_j) \mid \theta=d\psi \quad \mbox{for some}
\quad \psi\in W_{p,q}^{j-1}(M,\sigma_{j-1},\sigma_j)\}.
$$
The spaces
$$
H^j_{p,q}(M,\sigma_{j-1},\sigma_j)=Z_q^j(M,\sigma_j)/B^j_{p,q}(M,\sigma_{j-1},\sigma_j)
$$
and
$$
\overline{H}^j_{p,q}(M,\sigma_{j-1},\sigma_j)=Z_q^j(M,\sigma_j)
/\overline{B}^j_{p,q}(M,\sigma_{j-1},\sigma_j),
$$
where $\overline{B}^j_{p,q}(M,\sigma_{j-1},\sigma_j)$ is the closure of
$B^j_{p,q}(M,\sigma_{j-1},\sigma_j)$ in $L_q^j(M,\sigma_j)$ (equivalently,
in $W_{q,q}^j(M,\sigma_j,\sigma_j)$ are called
the $j$th {\it $L_{p,q}$-cohomology} and the $j$th {\it reduced $L_{p,q}$-cohomology}
of the Riemannian manifold~$M$ {\it with weights~$\sigma_{j-1}$ and~$\sigma_j$}.
The quotient space
$$
T^j_{p,q}(M,\sigma_{j-1},\sigma_j)=\overline{B}^j_{p,q}(M,\sigma_{j-1},\sigma_j)
/B^j_{p,q}(M,\sigma_{j-1},\sigma_j)
$$
will be referred to as the {\it $L_{p,q}$-torsion} of~$M$ with the given weights.
Clearly, the space $T^j_{p,q}(M,\sigma_{j-1},\sigma_j)$ is isomorphic to the closure
of the zero in $H^j_{p,q}(M,\sigma_{j-1},\sigma_j)$.  

Given a subset
$A\subset M$, the {\it relative nonreduced and reduced $L_{p,q}$-cohomology spaces }
$H^j_{p,q}(M,A,\sigma_{j-1},\sigma_j)$
and $\overline{H}^j_{p,q}(M,A,\sigma_{j-1},\sigma_j)$ are defined like
$$
H^j_{p,q}(M,A,\sigma_{j-1},\sigma_j)=Z_q^j(M,A,\sigma_j)
/B^j_{p,q}(M,A,\sigma_{j-1},\sigma_j)
$$
and
$$
\overline{H}^j_{p,q}(M,A,\sigma_{j-1},\sigma_j)=Z_q^j(M,A,\sigma_j)
/\overline{B}^j_{p,q}(M,A,\sigma_{j-1},\sigma_j),
$$
where the relative spaces $Z_q^j(M,A,\sigma_j)$ and $B^j_{p,q}(M,A,\sigma_{j-1},\sigma_j)$
are defined as their absolute analogs above with the spaces $W_{p,q}^j(M,\sigma_j,\sigma_{j+1})$
and $W_{p,q}^{j-1}(M,\sigma_{j-1},\sigma_j)$ replaced by 
$W_{p,q}^j(M,A,\sigma_j,\sigma_{j+1})$ and $W_{p,q}^{j-1}(M,A,\sigma_{j-1},\sigma_j)$.

For $p=q$, we write the subscript $p$ instead of~$p,p$ throughout. If the weights involved
in the definition of the corresponding space are equal to~$1$ then they will be omitted.

The spaces $W_{p,q}$ and $L_{p,q}$-cohomology were introduced at the beginning
of the 1980's by Gol$'$dshtein, Kuz$'$minov, and
Shvedov~\cite{GKSh821,GKSh822,GKSh83,GKSh84,GKSh85,GKSh89},
who obtained many results concerning $W_{p,q}$-forms and especially $L_p$-cohomology. Later
$L_{p,q}$-cohomology was considered in~\cite{GT98,GT06,Ko07}. 

In this paper, we, following~\cite{GKSh901,GKSh902}, look for conditions of the nontriviality
of the $L_{p,q}$-cohomology and $L_{p,q}$-torsion on warped cylinders, a class
of warped products of Riemannian manifolds.
By the {\it warped product} $X\times_f Y$ of two Riemannian manifolds $(X,g_X)$ and
$(Y,g_Y)$ {\it with the warping function}~$f:X\to\mathbb{R}_+$ we mean
the product manifold $X\times Y$ endowed with the metric $g_X+f^2(x)g_Y$. If $X=[a,b[$
is a half-interval on the real line then $X\times_f Y$ is referred to as
the {\it warped cylinder}. The study $L_2$-cohomology of warped cylinders was initiated by Cheeger~\cite{Ch}.

The structure of the article is as follows. In Section~\ref{half-interval}, we adapt
the results of~\cite{GKSh902} about the $L_p$-cohomology of a half-interval to the case $p\ne q$.
After that, using these $L_{p,q}$-results, in Section~\ref{Lpq_warped}, we prove a partial
$L_{p,q}$-generalization of Theorem~1 of~\cite{GKSh902} about the $L_p$-cohomology of
a~warped cylinder $[a,b[\times_f Y$ depending on the analytic properties of the function~$f$.
As an application, we obtain an extension of the necessary condition for the triviality of
the $L_{p,q}$-torsion of a surface of revolution in~$\mathbb R^{n+2}$~\cite{Ko99}
from the case $p=q$ to arbitrary $p$,$q$ such that
$\frac{1}{q}-\frac{1}{p}<\frac{1}{n+1}$.

\section{Weighted $L_{p,q}$-Cohomology of a Half-Interval}
\label{half-interval}
Consider a half-interval $[a,b[$, $-\infty< a<b \le\infty$ and positive continuous
functions $v_0, v_1:[a,b[\to\mathbb{R}$. For $1<{p,q}<\infty$, the space
$W_{p,q}^0([a,b[,v_0,v_1)$ can be identified with the space of the functions
$g\in L_p([a,b[, v_0)$ whose weak derivative $g'\in L_q([a,b[,v_1)$. As above, endow
$W_{p,q}^0([a,b[,v_0,v_1)$ with the norm
$$
\|g\|_{W_{p,q}^0([a,b[,v_0,v_1)}
= \biggl( \int_a^b |g(t)|^p v_0^p dt \biggr)^{1/p}
+ \biggl( \int_a^b |g'(t)|^q v_0^q dt \biggr)^{1/q}.
$$

From the classical Sobolev Embedding Theorem it follows that
the functions of the class $W_{p,q}^0([a,b[,v_0,v_1)$ are continuous
on $[a,b[$. Consider also the space
$$
W_{p,q}^0([a,b[,\{a\}, v_0,v_1)=\{f\in W_{p,q}^0([a,b[,\{a\},v_0,v_1\mid f(a)=0 \}.
$$
We have
\begin{gather*}
H_{p,q}^1([a,b[,v_0,v_1)=W_q^1([a,b[,v_1,v_1)/dW_{p,q}^0([a,b[,v_0,v_1); \\
H_{p,q}^1([a,b[,\{a\},v_0,v_1)=W_q^1([a,b[,\{a\},v_1,v_1)/dW_{p,q}^0([a,b[,\{a\},v_0,v_1).
\end{gather*}
The spaces $\overline{H}_{p,q}^1([a,b[,v_0,v_1)$ and
$\overline{H}_{p,q}^1([a,b[,\{a\},v_0,v_1)$ are described similarly.

We call the following assertion the lemma about the Hardy inequality~\cite{BS,GKSh901,M}:

\begin{lemma}\label{har}
Suppose that $1\le p,q \le\infty$, $\frac{1}{q}+\frac{1}{q'}=1$,
$\alpha,\beta\in [-\infty,\infty]$, $I_{\alpha,\beta}$ is the interval with endpoints
$\alpha$ and $\beta$, $v_0$ and $v_1$ are continuous positive functions on $I_{\alpha,\beta}$.
Then for the existence of a global constant~$C$ such that
$$
\biggl| \int_\alpha^\beta \biggl| v_0(t) \int_\alpha^\tau
g(t)dt \biggr|^p d\tau \biggr|^{1/p}
\le C \biggl| \int_\alpha^\beta |v_1(t)g(t)|^q dt \biggr|^{1/q},
$$
it is necessary and sufficient that $\chi_{p,q}(\alpha,\beta,v_0,v_1)<\infty$.

Here
$$
\chi_{p,q}(\alpha,\beta,v_0,v_1)
=\sup\limits_{\tau\in I_{\alpha,\beta}} \biggl\{ \biggl| \int_\tau^\beta
|v_0(t)|^p dt \biggr|^{1/p}
\biggl| \int_\alpha^\tau |v_1(t)|^{-q'}dt\biggr|^{1/q'} \biggr\}
$$
if $p\ge q${\rm;}
$$
\chi_{p,q}(\alpha,\beta,v_0,v_1)
= \Biggl| \int_\alpha^\beta \biggl( \biggl| \int_\alpha^\tau
|v_1(t)|^{-q'} dt\biggr|^{p-1} \biggl|\int_\tau^\beta |v_0(t)|^p dt \biggr|
\biggr)^{\frac{q}{q-p}} |v_1(\tau)|^{-q'} d\tau \biggr|^{\frac{q-p}{pq}}
$$
if $p<q$.

If $p=1$ {\rm($q'=\infty$)} then the corresponding integral 
must be replaced by $\mbox{\rm ess\,sup}$.
\end{lemma}

The constant $\chi_{p,q}(\alpha,\beta,v_0,v_1)$ will be referred to as the
\textit{Hardy constant}.

The following lemma was proved in~\cite{GKSh902} for $p=q$ and $v_0=v_1$. The proof given
in~\cite{GKSh902} holds
for different $p$ and $q$ and different $v_0$ and $v_1$.

\begin{lemma}\label{violhar}
Suppose that $\alpha,\beta\in [-\infty,\infty]$, $v_0,v_1:I_{\alpha,\beta}\to\mathbb{R}$
are positive continuous functions, and $\chi_{p,q}(\alpha,\beta,v_0,v_1)=\infty$. Then
there exists a nonnegative function $h$ such that
$$
\biggl|\int_\alpha^\beta v_1^q(t) h^q(t)dt|\biggr|<\infty, \quad
\biggl| \int_\alpha^\beta v_0^p(\tau)
\biggl| \int_\alpha^\tau h(t)dt \biggr|^p d\tau \biggr| =\infty.
$$
\end{lemma}

As in~\cite{GKSh902}, Lemma~2 yields the following assertion.

\begin{theorem}\label{interv}
If $v_0$, $v_1$ are positive continuous functions on $[a,b[$ and $1<p,q<\infty$ then

{\rm(1)} $H_{p,q}^1([a,b[,\{a\},v_0,v_1)=0 \Longleftrightarrow
\chi_{p,q}(a,b,v_0,v_1)<\infty${\rm;}

{\rm(2)} $H_{p,q}^1([a,b[,v_0,v_1)=0 \Longleftrightarrow \chi_{p,q}(a,b,v_0,v_1)<\infty
\; \mbox{or} \; \chi_{p,q}(b,a,v_0,v_1)<~\infty$.
\end{theorem}

\medskip

Let
\begin{equation}\label{compl}
0\to A \overset\varphi\to B \overset\psi\to C \to 0
\end{equation}
be an exact sequence of Banach complexes, i.e., complexes in the category of
Banach spaces and bounded linear operators. Sequence~(\ref{compl}) yields
an exact sequence of the cohomology spaces
$$
\dots \to H^{k-1}(C) \overset{\partial}\to H^k(A) \overset{\varphi^*}\to H^k(B)
\overset{\psi^*}\to H^k(C) \to \dots
$$
with continuous operators $\partial^*$, $\varphi^*$, $\psi^*$ and a semi-exact sequence
of the reduced cohomology spaces
\begin{equation}\label{cohseq}
\dots \to \overline H^{k-1}(C) \overset{\overline\partial}
\to \overline H^k(A) \overset{\overline\varphi^*}\to
\overline H^k(B) \overset{\overline\psi^*}\to \overline H^k(C) \to \dots
\end{equation}

Under certain conditions, sequence~(\ref{cohseq}) is exact at some terms
(see~\cite{GKSh901,KoK03,KSh99}). In particular,
Gol$'$dshtein, Kuz$'$minov, and Shvedov proved the following assertion
in~\cite{GKSh901}:

\begin{lemma}\label{exact}
If $H^{k+1}(A)=\overline H^{k+1}(A)$ and $\dim\varphi^*(H^k(A))<\infty$ then the sequence
$\overline H^k(A) \overset{\overline\varphi^*}\to \overline H^k(B) \overset{\overline\psi^*}
\to \overline H^k(C)
\overset{\overline\partial}\to \overline H^{k+1}(A)$ is exact.
\end{lemma}

As was explained in~\cite{GT06}, we can describe the $j$th weighted $L_{p,q}$-cohomology
of an $n$-dimensional Riemannian manifold~$M$ with given weights~$\sigma_{j-1}$ and $\sigma_j$
in terms of Banach complexes. To this end, consider an arbitrary sequence
$\pi=\{p_0,p_1,\dots,p_n\}\subset [1,\infty]$ with $p_{k-1}=p$ and $p_k=q$
and a sequence of positive continuous weights $\sigma=\{\sigma_k\}_{k=0}^n$
with the given~$\sigma_{j-1}$ and $\sigma_j$. Given a subset $A\subset M$, put
$$
W_\pi^k(M,A,\sigma)=W_{p_k,p_{k+1}}(M,A,\sigma_k,\sigma_{k+1}).
$$
Here we have assumed that $p_{n+1}=p_n$ and $\sigma_{n+1}=\sigma_n$.

Since the exterior differential is a bounded operator
$d^{k-1}:W_\pi^{k-1}(M,A,\sigma)\to W_\pi^k (M,A,\sigma)$,
we obtain a Banach complex  
\begin{equation}\label{bancom}
0 \to W_\pi^0(M,A,\sigma) \overset{d^0}\to  W_\pi^1(M,A,\sigma) \to \dots
\overset{d^{n-1}}\to W_\pi^n(M,A,\sigma) \to 0.
\end{equation}
By the $k$th {\it $L_\pi$-cohomology $H_\pi^k(M,A,\sigma)$} ({\it reduced $L_\pi$-cohomology
$\overline H_\pi^k(M,A,\sigma)$}) {\it of the Riemannian manifold~$M$ with respect
to~$A$ with weight~$\sigma$}
we mean the cohomology (reduced cohomology) of~(\ref{bancom}). Thus,
$H_\pi^k(M,A,\sigma)=H_{p_{k-1},p_k}^k(M,A,\sigma_{k-1},\sigma_k)$ and
$\overline H_\pi^k(M,A,\sigma)=\overline H_{p_{k-1},p_k}^k(M,A,\sigma_{k-1},\sigma_k)$
for all~$k$. In particular,
$$
H_\pi^j(M,A,\sigma)=H_{p,q}^j(M,A,\sigma_{j-1},\sigma_j), \quad
\overline H_\pi^j(M,A,\sigma)=\overline H_{p,q}^j(M,A,\sigma_{j-1},\sigma_j).
$$

Take $M=[a,b[$, $A=\{a\}$, $1<p,q<\infty$, $\pi=\{p,q\}$, and a pair of weights
$v=\{v_0,v_1\}$. We have the following exact sequence of Banach complexes:
$$
0\to W_\pi^*([a,b[,\{a\},v) \overset{\mathfrak j}\to W_\pi^*([a,b[,v)
\overset{\mathfrak i}\to H^*(\{a\})\to 0,
$$
where $H^*(\{a\})$ is the complex with the only nontrivial term $H^0(\{a\})=\mathbb{R}$.
Lemma~\ref{exact} yields the exact sequence
$$
\mathbb R=H^0(\{a\})\overset{\overline\partial}\to
\overline H_{p,q}^1([a,b[,\{a\},v_0,v_1) \overset{\overline{\mathfrak j}^*}\to
\overline H_{p,q}^1([a,b[,v_0,v_1).
$$

Thus, we infer the following assertion, proved for~$p=q$ in~\cite{GKSh902}. With what has been
said above, the proof of~\cite{GKSh902} extends to the case of $p\ne q$ without change.

\begin{theorem}\label{redcom}
If $v_0$, $v_1$ are positive continuous functions on~$[a,b[$, $1<p<\infty$, $1<q<\infty$,
$q'=\frac{q}{q-1}$,
$\frac{1}{q}+\frac{1}{q'}=1$ then

{\rm(1)} $\overline H_{p,q}^1([a,b[,v_0,v_1)=0${\rm;}

{\rm(2)} $\overline H_{p,q}^1([a,b[,\{a\},v_0,v_1)=0$ if and only if
$\int_a^b v_1^{-q'}(t)dt=\infty$ or $\int_a^b v_0^p (t)dt<\infty${\rm;}

{\rm(3)} If $\overline H_{p,q}([a,b[,\{a\},v_0,v_1)\ne 0$ then
$\overline\partial:\mathbb{R}=H^0(\{a\}) \to H_{p,q}^1([a,b[,\{a\},v_0,v_1)$
is an isomorphism.
\end{theorem}

\section{$L_{p,q}$-Cohomology of the Warped Cylinder $C_{a,b}^f$}\label{Lpq_warped}

Let $Y$ be an orientable manifold of dimension~$n$,
$C_{a,b}^f Y=[a,b[\times_f Y$. Put $Y_a=\{a\}\times Y$. Generally speaking,
$C_{a,b}^f$ is a~Lipschitz Riemannian manifold in the sense of~\cite{GKSh821} but
we will assume for simplicity that $\partial Y=\varnothing$ to make $C_{a,b}^f$
smooth, which will be enough for our purposes.

Suppose that $1<p<\infty$ and $1<q<\infty$.

In~\cite{GKSh902}, Gol$'$dshtein, Kuz$'$minov, and Shvedov introduced the bilinear mapping
$$
\nu: L_p^{j-1}(Y) \times L_p^1([a,b[,f^{\frac{n}{p}-j+1})\to L_p^j(C_{a,b}^f Y),
$$
$\nu(\varphi,gdt)=gdt\wedge\varphi$. In~\cite{GKSh902} it was proved that $\nu$ is
continuous and if $\varphi\in Z_p^{j-1}(Y)$ then
$\nu_\varphi=\nu(\varphi,\cdot):L_p^1([a,b[,f^{\frac{n}{p}-j+1})\to L_p^j(C_{a,b}^f Y)$
induces continuous mappings
\begin{gather*}
\nu_\varphi^*:H_p^1([a,b[,f^{\frac{n}{p}-j+1}) \to H_p^j(C_{a,b}^f Y); \\
\tilde\nu_\varphi^*:H_p^1([a,b[,\{a\},f^{\frac{n}{p}-j+1})
\to H_p^j(C_{a,b}^f Y,Y_a).
\end{gather*}

Supposing that $\varphi\in Z_p^{j-1}(Y)\cap Z_q^{j-1}(Y)$, we similarly become
convinced that the mapping $\nu_\varphi=\nu(\varphi,\cdot)$ induces continuous mappings
\begin{gather*}
\nu_\varphi^*:H_{p,q}^1([a,b[,f^{\frac{n}{p}-j+1},f^{\frac{n}{q}-j+1})
\to H_{p,q}^j(C_{a,b}^f Y); \\
\tilde\nu_\varphi^*:H_{p,q}^1([a,b[,\{a\},f^{\frac{n}{p}-j+1},f^{\frac{n}{q}-j+1})
\to H_{p,q}^j(C_{a,b}^f Y,Y_a).
\end{gather*}

Now, assume that $\psi\in L_{p'}^{n+1-j}(Y)$ ($p'=\frac{p}{p-1}$) and
$\omega\in L_p^j(C_{a,b}Y)$. Write $\omega$ in the form
$\omega=\omega_A+dt\wedge\omega_B$, where $\omega_A$, $\omega_B$ do not
contain~$dt$~\cite{GKSh901}. Following~\cite{GKSh902}, introduce the continuous operator
$\mu_\psi:L_p^j(C_{a,b}^f Y) \to L_p^1([a,b[,f^{\frac{n}{p}-j+1})$ by the formula
$$
\mu_\psi\omega = \biggl(\int_Y \omega_B(t)\wedge\psi \biggr) dt.
$$

The following lemma was proved in~\cite{GKSh902} for $p=q$ and
$\psi\in V_{p'}^{n-j+1}(Y)$. The proof in~\cite{GKSh902} easily extends
to~$p\ne q$:

\begin{lemma}
If $\psi\in D^{n-j+1}(Y)$ and $d\psi=0$ then $\mu_\psi$ induces continuous mappings
\begin{gather*}
\mu_\psi^*: H_{p,q}^j(C_{a,b}^f Y) \to
H_{p,q}^1 ([a,b[,f^{\frac{n}{p}-j+1},f^{\frac{n}{q}-j+1}); \\
\tilde\mu_\psi^*: H_{p,q}^j(C_{a,b}^f Y, Y_a) \to
H_{p,q}^1 ([a,b[,\{a\},f^{\frac{n}{p}-j+1},f^{\frac{n}{q}-j+1})
\end{gather*}
\end{lemma}

We have the following theorem partially generalizing item~7 of Theorem~1 in~\cite{GKSh902}:

\begin{theorem}\label{main}
Suppose that $Y$ is an orientable $n$-dimensional Riemannian manifold,
$\infty<a<b\le\infty$, $f:[a,b[\to\mathbb{R}$ is a positive continuous function,
$1<p<\infty$, $1<q<\infty$. Assume that there exists
$\varphi\in Z_p^{j-1}(Y)\cap Z_q^{j-1}(Y)$ such that
$\int_Y \varphi\wedge\gamma\ne 0$ for some form $\gamma\in D^{n-j+1}(Y)$,
$d\gamma=0$.

The following hold{\rm:}

{\rm(1)} if $\chi_{p,q}(a,b,f^{\frac{n}{p}-j+1},f^{\frac{n}{q}-j+1})=\infty$ then
$H_{p,q}^j(C_{a,b}^f Y,Y_a)\ne 0${\rm;}

{\rm(2)} if $\chi_{p,q}(a,b,f^{\frac{n}{p}-j+1},f^{\frac{n}{q}-j+1})=\infty$
and $\chi_{p,q}(b,a,f^{\frac{n}{p}-j+1},f^{\frac{n}{q}-j+1})=\infty$ then
$T_{p,q}^j(C_{a,b}^f)\ne 0$ and, hence, $\dim H_{p,q}^j(C_{a,b}^f Y)=\infty$.
\end{theorem}

\begin{proof}
Let $\varphi\in Z_p^{j-1}(Y)\cap Z_q^{j-1}(Y)$ be a cocycle having the property mentioned
in the theorem and let $\gamma\in D^{n-j+1}(M)$ be a form such that
$\int_Y \varphi\wedge\gamma=1$. Then
$\mu_\gamma^*\circ \nu_\varphi^* =\id$,
$\tilde\mu_\gamma^*\circ\tilde\nu_\varphi^*=\id$~\cite{GKSh902}. Consequently,
the mappings
$\nu_\varphi^*:H_{p,q}^1([a,b[,f^{\frac{n}{p}-j+1},f^{\frac{n}{q}-j+1})
\to H_{p,q}^1(C_{a,b}^f Y)$
and
$\tilde\nu_\varphi^*:H_{p,q}^1([a,b[,\{a\},f^{\frac{n}{p}-j+1},f^{\frac{n}{q}-j+1})
\to H_{p,q}^1(C_{a,b}^f Y,Y_a)$
are injective.

Suppose that $\chi_{p,q}(a,b,f^{\frac{n}{p}-j+1},f^{\frac{n}{q}-j+1})=\infty$.
Then, by Theorem~\ref{interv}, \linebreak
$H_{p,q}^1([a,b[,\{a\},f^{\frac{n}{p}-j+1},f^{\frac{n}{q}-j+1})\ne 0$. Therefore,
$H_{p,q}^j(C_{a,b}^f,Y_a)\ne 0$.

Assume now that $\chi_{p,q}(a,b,f^{\frac{n}{p}-k+1},f^{\frac{n}{q}-j+1})=\infty$
and $\chi_{p,q}(b,a,f^{\frac{n}{p}-k+1},f^{\frac{n}{q}-j+1})=\infty$. Then,
by Theorem~\ref{interv},
$H_{p,q}^1([a,b[,f^{\frac{n}{p}-j+1},f^{\frac{n}{q}-j+1})\ne 0$. Since,
by Theorem~\ref{redcom},
$\overline H_{p,q}([a,b[,f^{\frac{n}{p}-j+1},f^{\frac{n}{q}-j+1})= 0$, we have
$T_{p,q}^j([a,b[,f^{\frac{n}{p}-j+1},f^{\frac{n}{q}-j+1})\ne 0$. Hence,
$T_{p,q}^j(C_{a,b}^f Y)\ne 0$. The theorem is proved.
\end{proof}

\section*{$L_{p,q}$-Torsion of a Surface of Revolution}

Let $M$ be a surface of revolution in~$\mathbb{R}^{n+2}$, i.e.,
the ($n+1$)-dimensional surface defined by the equation
\begin{equation}\label{surfrev}
f^2(x_1)=x_2^2+\dots+x_{n+2}^2, \quad (x_1,\dots,x_{n+2})\in\mathbb{R}^{n+2}, \,\,
x_1\ge 0,
\end{equation}
where $f:[0,\infty[\to\mathbb{R}$ is a positive smooth function. With the metric induced
from~$\mathbb{R}^{n+2}$, the manifold~$M$  is the product $[0,\infty[\times \mathbb{S}^n$
with the metric
$$
g_M=(1+f'^2(x_1))dx_1^2+f^2(x_1)dy^2,
$$
where $dx_1^2$ and $dy^2$ are the conventional Riemannian metrics on~$[0,\infty[$ and
the sphere~$\mathbb{S}^n$. In other words, $M$ may be considered as the warped product
$[0,\infty[\times_F \mathbb{S}^n$, where $F=f\circ G^{-1}$,
$G(x)=\int_0^x \sqrt{1+f'^2(t)}dt$.

In~\cite{Ko07}, we have proved the following fact:

\begin{theorem}\label{normsol1}
Suppose that $f$ is unbounded, $p,q\in [1,\infty[$, $\frac{1}{q}-\frac{1}{p}<\frac{1}{n+1}$,
$1\le j \le n+1$. Then $T_{p,q}^{j}(M)\ne 0$.
\end{theorem}

Kuz$'$minov and Shvedov~\cite{KSh96} established that $T_p^j(M)$ is zero
for all $j$, $2\le j \le n$ and that, for~$j=1,n+1$, the triviality
of $T_p^j(M)$ depends on the finiteness of some Hardy constants. This is due
to the connection between the $L_p$-cohomology of the warped product $C_{a,b}^f Y$
and the weighted $L_p$-cohomology of $[a,b[$ given in the mentioned
papers~\cite{GKSh901,GKSh902}.
Above we have shown that there is a connection of this type for $L_{p,q}$-cohomology.
Namely, by Theorem~\ref{main}, since $\mathbb{S}^n$ is compact and the de Rham
cohomology~$H^{j-1}(\mathbb{S}^n)$ of~$\mathbb{S}^n$ is nontrivial if~$j=1,n+1$,
for $T_{p,q}^j(M)$ ($j=1,n+1$) to be zero, it is necessary that
$\chi_{p,q}(0,\infty,F^{\frac{n}{p}-j+1},F^{\frac{n}{q}-j+1})<\infty$
or
$\chi_{p,q}(\infty,0,F^{\frac{n}{p}-j+1},F^{\frac{n}{q}-j+1})<\infty$.

The main result of this section is a generalization of Theorems~2 and~2$'$ of~\cite{Ko99}
and is formulated as follows:

\begin{theorem}\label{tors}
Let $M$ be the surface of revolution~{\rm(\ref{surfrev})}. Suppose that $1<p<\infty$,
$1<q<\infty$,
$\frac{1}{q}-\frac{1}{p}<\frac{1}{n+1}, j\in \{1,n+1\}$. If $T_{p,q}^j(M)=0$ then
$\lim\limits_{x\to\infty} f(x)=0$ and $\vol M<\infty$.
\end{theorem}

\begin{proof}

Put $k=j-1$.

We have the following equalities:
\begin{multline*}
\chi_{p,q}^0\equiv \chi_{p,q}(0,\infty,F^{\frac{n}{p}-k},F^{\frac{n}{q}-k})
\\
=\sup\limits_{\tau>0} \biggl\{ \biggl( \int_\tau^\infty
f^{n-kp}(t) \sqrt{1+f'^2(t)} dt \biggr)^{1/p}
\biggl( \int_0^\tau f^{-(\frac{n}{q}-k)q'}
\sqrt{1+f'^2(t)} dt \biggr)^{1/q'} \biggr\};
\end{multline*}
\begin{multline*}
\chi_{p,q}^\infty\equiv \chi_{p,q}(\infty,0,F^{\frac{n}{p}-k},F^{\frac{n}{q}-k})
\\
=\sup\limits_{\tau>0} \biggl\{ \biggl( \int_0^\tau
f^{n-kp}(t) \sqrt{1+f'^2(t)} dt \biggr)^{1/p}
\biggl( \int_\tau^\infty f^{-(\frac{n}{q}-k)q'}
\sqrt{1+f'^2(t)} dt \biggr)^{1/q'} \biggr\}
\end{multline*}
if $p\ge q$;
\begin{multline}\label{int1}
\chi_{p,q}^0\equiv \chi_{p,q}(0,\infty,F^{\frac{n}{p}-k},F^{\frac{n}{q}-k})
\\
= \Biggl( \int_0^\infty \biggl[ \biggl( \int_0^{H(x)}
f^{-(\frac{n}{q}-k)q'}(t) \sqrt{1+f'^2(t)} dt \biggr)^{p-1}
\int_{H(x)}^\infty f^{n-kp}(t)\sqrt{1+f'^2(t)}dt \biggr]^{\frac{q}{q-p}}
\\
\times f^{-(\frac{n}{q}-k)q'}(x)\sqrt{1+f'^2(x)}dx \Biggr)^{\frac{q-p}{qp}};
\end{multline}
\begin{multline*}
\chi_{p,q}^\infty \equiv \chi_{p,q}(\infty,0,F^{\frac{n}{p}-k},F^{\frac{n}{q}-k})
\\
= \Biggl( \int_0^\infty \biggl[ \biggl( \int_{H(x)}^\infty
f^{-(\frac{n}{q}-k)q'}(t) \sqrt{1+f'^2(t)} dt \biggr)^{p-1}
\int_0^{H(x)} f^{n-kp}(t)\sqrt{1+f'^2(t)}dt \biggr]^{\frac{q}{q-p}}
\\
\times f^{-(\frac{n}{q}-k)q'}(x)\sqrt{1+f'^2(x)}dx \Biggr)^{\frac{q-p}{qp}}
\end{multline*}
if $p<q$.
Here $H(x)$ is the function inverse to the arc length function
$G(x)=\int_0^x \sqrt{1+f'^2(t)}dt$.

The main element in the proof of Theorem~\ref{tors} is the following lemma which
has some independent interest.

\begin{lemma}\label{chipq}
If $\frac{1}{q}-\frac{1}{p}<\frac{1}{n+1}$, $1<p<\infty$, $1<q<\infty$,
$0\le k\le n$, then the following hold{\rm:}

{\rm(1)} if $\chi_{p,q}^0<\infty$ or $\chi_{p,q}^\infty<\infty$ then
$\lim\limits_{t\to\infty}f(t)=0${\rm;}

{\rm(2)} if $\frac{n}{p}-k\le 0$ then $\chi_{p,q}^0=\infty${\rm;}

{\rm(3)} if $\frac{n}{q}-k\ge 0$ then $\chi_{p,q}^\infty=\infty$.
\end{lemma}

\begin{proof}
Suppose first that $p\ge q$.

Assume that $\chi_{p,q}^0<\infty$. Then
\begin{equation}\label{fin}
\int_0^\infty f^{n-kp}(t)\sqrt{1+f'^2(t)}dt<\infty
\end{equation}
Since
$$
f^{n-kp}(t)\sqrt{1+f'^2(t)}\ge f^{n-kp}(t)|f'(t)|,
$$
it follows that the integral
\begin{multline}\label{fin1}
\int_0^\infty f^{n-kp}(t)f^\prime(t)dt
\\
= \begin{cases}
\frac{1}{n-kp+1} \lim\limits_{t\to\infty}
(f^{n-kp+1}(t)-f^{n-kp+1}(0)) &\text{if $n-kp\ne -1$}, \\
\lim\limits_{t\to\infty} \log\frac{f(t)}{f(0)}
&\text{if $n-kp=-1$}
\end{cases}
\end{multline}
is finite.

There appear several possibilities:

(a) $\frac{n}{p}-k>0$.
The above implies that there exists a finite limit
$\lim\limits_{t\to\infty}f(t)$,
which is zero by~(\ref{fin}).

(b) $\frac{n}{p}-k=0$.
This is impossible in view of~(\ref{fin}).

(c) $-\frac{1}{p}<\frac{n}{p}-k<0$.
Then
$n-kp+1>0$
and
$f(t)$
has a finite limit as
$t\to\infty$,
which contradicts~(\ref{fin}).

(d) $\frac{n}{p}-k=-\frac{1}{p}$.
A contradiction to~(\ref{fin}).

(e) $\frac{n}{p}-k<-\frac{1}{p}$.
In this case,
$n-kp<-1$.
Hence,
$\lim\limits_{t\to\infty}f(t)=\infty$.

Note that, since $\frac{1}{q}-\frac{1}{p}<\frac{1}{n+1}$, we have
$k+1>\frac{n+1}{p}+1>\frac{n+1}{q}$, whence
$-(\frac{n}{q}-k)q'+1>0$. We infer
\begin{multline}\label{chain2}
\biggl( \int_\tau^\infty
f^{n-kp}(t) \sqrt{1+f'^2(t)} dt \biggr)^{1/p}
\biggl( \int_0^\tau f^{-(\frac{n}{q}-k)q'}
\sqrt{1+f'^2(t)} dt \biggr)^{1/q'}
\\
\ge
\biggl( \int_\tau^\infty
f^{n-kp}(t) |f'(t)| dt \biggr)^{1/p}
\biggl( \int_0^\tau f^{-(\frac{n}{q}-k)q'}
|f'(t)| dt \biggr)^{1/q'}
\\
\ge
\biggl| \int_\tau^\infty
f^{n-kp}(t) f'(t) dt \biggr|^{1/p}
\biggl| \int_0^\tau f^{-(\frac{n}{q}-k)q'}
f'(t) dt \biggr|^{1/q'}
\\
\ge
\biggl( \frac{f^{n-kp+1}(\tau)}{|n-kp+1|}\biggr)^{1/p}
\biggl| \frac{f^{-\frac{n}{q}-k)q'+1}(\tau)-f^{-(\frac{n}{q}-k)q'+1}(0)}
{-(\frac{n}{q}-k)q'+1} \biggr|^{1/q'}
\\
= C\cdot f^{\frac{n+1}{p}-\frac{n+1}{q}+1}(\tau)
|1-f^{-(\frac{n}{q}-k)q'+1}(0)f^{(\frac{n}{q}-k)q'-1}(\tau)|^{1/q'}.
\end{multline}
The last quantity in~(\ref{chain2}) is equivalent to
$Cf^{\frac{n+1}{p}-\frac{n+1}{q}+1}(\tau)$ as $\tau\to\infty$ and, hence,
tends to infinity. Therefore, $\chi_{p,q}^0=\infty$, and we obtain a contradiction.

Thus, if $\chi_{p,q}^0<\infty$ then $\lim\limits_{t\to 0}f(t)=0$ and $\frac{n}{p}-k>0$.

Suppose now that $\chi_{p,q}^\infty<\infty$. Then
\begin{equation}\label{fin2}
\int_0^\infty f^{-(\frac{n}{q}-k)q'}(t)\sqrt{1+f'^2(t)}dt<\infty
\end{equation}
and, hence, there exists a finite integral
\begin{multline}\label{fin3}
\int_0^\infty f^{-(\frac{n}{q}-k)q'}(t)f^\prime(t)dt
\\
= \begin{cases}
\frac{\lim\limits_{t\to\infty}
f^{-(\frac{n}{q}-k)q'+1}(t)-f^{-(\frac{n}{q}-k)q'+1}(0)}{-(\frac{n}{q}-k)q'+1}
&\text{if $-(\frac{n}{q}-k)q'\ne -1$}, \\
\lim\limits_{t\to\infty} \log\frac{f(t)}{f(0)}
&\text{if $-(\frac{n}{q}-k)q'= -1$}.
\end{cases}
\end{multline}
As in the case $\chi_{p,q}^0<\infty$, we infer that either $\frac{n}{q}-k<0$ and
$\lim\limits_{t\to\infty}f(t)=0$ or $(\frac{n}{q}-k)q'>1$ and
$\lim\limits_{t\to\infty}f(t)=\infty$. In the latter case we have:
\begin{multline}\label{chain3}
\chi_{p,q}^\infty \ge \sup\limits_{\tau>0} \biggl\{
\biggl( \frac{f^{-(\frac{n}{q}-k)q'+1}(\tau)}{1-(\frac{n}{q}-k)q'}\biggr)^{1/q'}
\biggl| \frac{f^{n-kp+1}(\tau)-f^{n-kp+1}(0)}{n-kp+1} \biggr|^{1/p} \biggl\}
\\
= C\, \sup\limits_{\tau\to 0} \bigl\{ f^{\frac{n+1}{p}-\frac{n+1}{q}+1}(\tau)
|1-f^{n-kp+1}(0) f^{-(n-kp+1)}(\tau)|^{1/p}  \bigr\},
\end{multline}
where $C=\const>0$. Since $k<\frac{n+1}{q}-1<\frac{n+1}{p}$, we have $n-kp+1>0$,
and, hence, the last quantity in~(\ref{chain3}) behaves like
$Cf^{\frac{n+1}{p}-\frac{n+1}{q}+1}(\tau)$ and, consequently, tends to infinity
as $\tau\to\infty$. Hence, $\chi_{p,q}^\infty=\infty$; a contradiction.

Thus, if $\chi_{p,q}^\infty<\infty$ then $\lim\limits_{t\to 0}f(t)=0$ and $\frac{n}{q}-k<0$.

We now pass to the case $p<q$.

Suppose that $\chi_{p,q}^0<\infty$. Then, as above, we have~(\ref{fin}) and~(\ref{fin1})
and conclude that either $\frac{n}{p}-k>0$ and $\lim\limits_{t\to\infty}=0$
or $\frac{n}{p}-k<-\frac{1}{p}$ and $\lim\limits_{t\to\infty}f(t)=\infty$. Show that
the latter case is impossible. By~(\ref{int1}), we infer
\begin{multline}\label{chipq0}
(\chi_{p,q}^0)^{\frac{pq}{q-p}} \ge
\int_0^\infty \biggl[ \biggl| \int_0^{H(x)}
f^{-(\frac{n}{q}-k)q'}(t) f'(t) dt \biggr|^{p-1}
\biggl|\int_{H(x)}^\infty f^{n-kp}(t) f'(t) dt\biggr| \;
\biggr]^{\frac{q}{q-p}}
\\
\times f^{-(\frac{n}{q}-k)q'}(x)\sqrt{1+f'^2(x)}dx
\\
=\int_0^\infty \biggl[ \biggl|
\frac{f^{-(\frac{n}{q}-k)q'+1}(H(x))-f^{-(\frac{n}{q}-k)q'+1}(0)}{-(\frac{n}{q}-k)q'+1}
\biggr|^{p-1}
\biggl| \frac{f^{n-kp+1}(H(x))}{n-kp+1}\biggr|\biggr]^{\frac{q}{q-p}}
\\
\times f^{-(\frac{n}{q}-k)q'}(x) \sqrt{1+f'^2(x)} dx
\\
=\int_0^\infty  \biggl[ \biggl|
\frac{F^{-(\frac{n}{q}-k)q'+1}(s)-F^{-(\frac{n}{q}-k)q'+1}(0)}{-(\frac{n}{q}-k)q'+1}
\biggr|^{p-1}
\biggl| \frac{F^{n-kp+1}(s)}{n-kp+1}\biggr|\biggr]^{\frac{q}{q-p}}
F^{-(\frac{n}{q}-k)q'}(s) ds
\\
= C \int_0^\infty F^N(s)
|1-F^{-(\frac{n}{q}-k)q'+1}(0)F^{(\frac{n}{q}-k)q'-1}(s)| ds.
\end{multline}

Here $C=\const>0$ and
\begin{multline*}
N=\biggl(\biggl(-\biggl(\frac{n}{q}-k\biggr)q'+1\biggr)(p-1)+n-kp+1\biggr) \frac{q}{q-p}
-\biggl(\frac{n}{q}-k\biggr)q'
\\
=\biggl[ \biggl( 1-\frac{p-1}{q-1} \biggr) \frac{q}{q-p} - \frac{1}{q-1} \biggr] n
- \biggl[ \biggl( \frac{q(p-1)}{q-1}-p \biggr) \frac{q}{q-p} + \frac{1}{q-1} \biggr] k
+\frac{pq}{q-p}
\\
= n - k + \frac{pq}{q-p} > 0.
\end{multline*}

Moreover, $(\frac{n}{q}-k)q'-1<0$, since  $\frac{n}{q}<\frac{n+1}{q}<\frac{n+1}{p}<k$.
Consequently, the expression under the last integral in~(\ref{chipq0}) is equivalent
to $CF^{n-k+\frac{pq}{q-p}}(s)$, i.e., tends to~$\infty$ as $s\to\infty$
and, thus, the integral does not exist. A contradiction.

Suppose now that $\chi_{p,q}^\infty<\infty$. Then we have~(\ref{fin2}) and~(\ref{fin3})
and infer that, in this case, either $\frac{n}{q}-k<0$ and
$\lim\limits_{t\to\infty} f(t)=0$ or $q'(\frac{n}{q}-k)>1$ and
$\lim\limits_{t\to\infty}f(t)=\infty$. In the latter case, we infer
\begin{multline}\label{chipqin}
(\chi_{p,q}^\infty)^{\frac{pq}{q-p}} \ge
\int_0^\infty \biggl[ \biggl| \int_0^{H(x)}
f^{n-kp}(t) f'(t) dt \biggr| \;
\biggl|\int_{H(x)}^\infty f^{-(\frac{n}{q}-k)q'}(t) f'(t) dt \biggr|^{p-1} \;
\biggr]^{\frac{q}{q-p}}
\\
\times f^{-(\frac{n}{q}-k)q'}(x)\sqrt{1+f'^2(x)}dx
\\
=\int_0^\infty \biggl[ \biggl|
\frac{f^{n-kp+1}(H(x))-f^{n-kp+1}(0)}{n-kp+1}
\biggr|
\biggl| \frac{f^{-(\frac{n}{q}-k)q'+1}(H(x))}{-(\frac{n}{q}-k)q'+1}\biggr|^{p-1}
\biggr]^{\frac{q}{q-p}}
\\
\times f^{-(\frac{n}{q}-k)q'}(x) \sqrt{1+f'^2(x)} dx
\\
=\int_0^\infty  \biggl[ \biggl|
\frac{F^{n-kp+1}(s)-F^{n-kp+1}(0)}{n-kp+1} \biggr|
\biggl| \frac{F^{-(\frac{n}{q}-k)q'+1}(s)}{-(\frac{n}{q}-k)q'+1}\biggr|^{p-1}
\biggr]^{\frac{q}{q-p}} F^{-(\frac{n}{q}-k)q'}(s) ds
\\
= C \int_0^\infty F^N(s)
|1-F^{n-kp+1}(0)F^{-(n-kp+1)}(s)| ds.
\end{multline}

Here, as above, $C=\const>0$, $N=n-k+\frac{pq}{q-p}>0$, and $-(n-kp+1)<0$.
Thus, the expression under the integral is equivalent to $CF^{n-k+\frac{pq}{q-p}}(s)$, i.e.,
tends to infinity as $s\to\infty$.

Lemma~\ref{chipq} is completely proved.
\end{proof}

Now, return to the proof of Theorem~\ref{tors}. Suppose that $T_{p,q}^j(M)=0$ for $j=1$
or $j=n+1$. Then, by Theorem~\ref{main},
$\chi_{p,q}(0,\infty,f^{\frac{n}{p}-j+1},f^{\frac{n}{q}-j+1})<\infty$ (and, hence,
$\int_0^\infty f^{(\frac{n}{p}-j+1)p(t)} \sqrt{1+f'^2(t)}dt<\infty$)
or $\chi_{p,q}(\infty,0,f^{\frac{n}{p}-j+1},f^{\frac{n}{q}-j+1})<\infty$
(and, hence, $\int_0^\infty f^{-(\frac{n}{q}-j+1)q'}(t) \sqrt{1+f'^2(t)}dt<\infty$).
By Lemma~\ref{chipq}, this implies that $\lim\limits_{t\to\infty}f(t)=0$ and, in both cases,
$$
\vol M = s_n\int_0^\infty f^n(t) \sqrt{1+f'^2(t)}dt<\infty.
$$
Here $s_n$ stands for the volume of the $n$-dimensional unit sphere in~$\mathbb R^{n+1}$.

The theorem is proved.
\end{proof}
\smallskip

\textbf{Acknowledgments.} The paper was begun in January 2007 during the author's research stay
at the Institut des Hautes \'Etudes Scientifiques (IH\'ES) in Bures-sur-Yvette. The author
would like to thank the Director Professor Jean-Pierre Bourguignon
for the invitation to visit the IH\'ES and all people at the Institute for their
warm hospitality. The author also expresses his gratitude to Vladimir Gol$'$dshtein,
Vladimir Kuz$'$\-mi\-nov, Aleksandr Romanov, and Igor Shvedov for useful discussions.

\end{document}